\documentclass[12pt]{article}
\usepackage{amsmath,amssymb,amsthm}
\usepackage{hyperref}
\usepackage{datetime}
\usepackage{units}
\usepackage{color}
\usepackage[T1]{fontenc}
\usepackage[utf8]{inputenc}
\usepackage{authblk}
\usepackage{bm,latexsym,mathrsfs,enumerate}
\usepackage{color}

\title{ The higher derivatives of the inverse tangent function and rapidly convergent BBP-type formulas\thanks{%
MSC 2010: 30D10, 40A25}}
\author[1]{Kunle Adegoke\thanks{adegoke00@gmail.com}}
\affil{Department of Physics and Engineering Physics, \mbox{Obafemi Awolowo University, Ile-Ife, 220005 Nigeria}}
\author[2]{Olawanle Layeni\thanks{olawanle.layeni@gmail.com }}
\affil{Department of Mathematics, \mbox{Obafemi Awolowo University, Ile-Ife, 220005 Nigeria}}

\theoremstyle{plain}
\numberwithin{equation}{section}

\begin{document}
\date{}
\maketitle
\begin{abstract}
\noindent We give a closed formula for the $n^{th}$ derivative of $\arctan x$. A new expansion for $\arctan x$ is also obtained and rapidly convergent series for $\pi$ and $\pi\sqrt 3$ are derived.
 
\end{abstract}
\tableofcontents

\section{Introduction}

The derivation of the $n^{th}$ derivative of $\arctan x$ is not straightforward (see e.g~\cite{euler93,forumarctan}). Efficient numerical computation of $\pi$ and related constants is often dependent on finding rapidly converging series for the inverse tangent function. Numerous interesting series and identities for $\pi$, ranging from the Gregory-Leibniz~\cite{eymard04} formula through the Machin-Like formulas~\cite{mathworldmachin} to the more recent BBP-type formulas~\cite{lord99}, are derived by manipulating the inverse tangent function. Of course there are also interesting series for $\pi$ whose connections with the arctangent function may not be obvious. Examples would be the numerous series for $\pi$, discovered by Ramanujan~\cite{berndt94}.

\smallskip

In this paper we will give a closed formula for the \mbox{$n^{th}$} derivative of $\arctan x$. A new series expansion for $\arctan x$ will also be obtained and rapidly convergent series for $\pi$ and $\pi\sqrt 3$ will be derived.

\section{The $n^{th}$ derivative of $\arctan x$ }\label{sec.derivative}
We have the following result.

THEOREM 1. The function $f(x)=\arctan x$ possesses derivatives of all order for \mbox{$x \in ( - \infty ,\infty )$}. The $n^{th}$ derivative of $\arctan x$ is given by the formula

\begin{equation}\label{equ.main}
\frac{{{\rm d}^n }}{{{\rm d}x^n }}\left( {\arctan x} \right) = \frac{{( - 1)^{n - 1} (n - 1)!}}{{(1 + x^2 )^{n/2} }}\sin \left[ {n\arcsin \left( {\frac{1}{\sqrt{1 + x^2 }}} \right)} \right],\quad n\in\mathbb{Z^+}\,.
\end{equation}

PROOF. It is convenient to make the right hand side of~\eqref{equ.main} more compact by writing
\[
\sin \theta  = \frac{1}{{\sqrt {1 + x^2 } }}\,.
\]
The formula then becomes
\[
\frac{{{\rm d}^n }}{{{\rm d}x^n }}\left( {\arctan x} \right) = ( - 1)^{n - 1} (n - 1)!\sin ^n \theta \sin n\theta\,. 
\]
The existence of the derivatives follows from the analyticity of $\arctan x$ on the real line. The proof of formula~\eqref{equ.main} is by mathematical induction. Clearly, the theorem is true for $n=1$. Suppose the theorem is true for $n=k$; that is, suppose

\begin{equation}\label{induction}
P_k:\quad\frac{{{\rm d}^k }}{{{\rm d}x^k }}\left( {\arctan x} \right) = ( - 1)^{k - 1} (k - 1)!\sin ^k \theta \sin k\theta\,. 
\end{equation}

We will show that the theorem is true for $n=k+1$ whenever it is true for $n=k$.

\bigskip

Differentiating both sides of~\eqref{induction} with respect to $x$ and noting that \mbox{${\rm d}\theta /{\rm d}x =  - \sin ^2 \theta $}, we have
\[
\begin{split}
\frac{{\rm d}}{{{\rm d}x}}\left[ {\frac{{{\rm d}^k }}{{{\rm d}x^k }}\left( {\arctan x} \right)} \right]& = ( - 1)^k k!\sin ^{k + 1} \theta \left( {\cos \theta \sin k\theta  + \cos k\theta \sin \theta } \right)\\
&= ( - 1)^k k!\sin ^{k + 1} \theta \sin (k + 1)\theta\,,
\end{split}
\]
so that
\[
P_{k+1}:\quad\frac{{{\rm d}^{k + 1} }}{{{\rm d}x^{k + 1} }}\left( {\arctan x} \right) = ( - 1)^k k!\sin ^{k + 1} \theta \sin (k + 1)\theta\,. 
\]
Thus $P_k  \Rightarrow P_{k + 1} $ and the proof is complete.

\section{A new expansion for $\arctan x$}
Perhaps the most well known series for $\arctan x$ is its Maclaurin expansion
\begin{equation}\label{equ.awqxz1m}
\begin{split}
\arctan x& = \sum_{n = 0}^\infty  {\frac{{( - 1)^n x^{2n + 1} }}{{2n + 1}}}\\
&= x - \frac{{x^3 }}{3} + \frac{{x^5 }}{5} - \frac{{x^7 }}{7} +  -  \cdots 
\end{split}
\end{equation}

Apart from its simplicity and elegance, series~\eqref{equ.awqxz1m} as it stands has little computational value as its radius of convergence is small~($R=1$) and the convergence is slow (logarithmic convergence) at the interesting endpoint~\mbox{$x=1$}. Note, however, that for $|x|<1$, one finds roughly $n\log_{10}(x)$ decimals of $\arctan x$, so that the convergence is linear. Euler transformation gives the form~\cite{castellanos88}:
\begin{equation}\label{equ.v0zvv88}
\arctan x = \sum_{n = 0}^\infty  {\frac{{2^{2n} (n!)^2 }}{{(2n + 1)!}}\frac{{x^{2n + 1} }}{{(1 + x^2 )^{n + 1} }}}\,. 
\end{equation}

The ratio test establishes easily that the series~\eqref{equ.v0zvv88} converges for all real $x$, giving $R=\infty$. Formula \eqref{equ.v0zvv88} exhibits linear convergence.

\bigskip

We now present a new series for $\arctan x$.

\smallskip

THEOREM 2. The function $f(x)=\arctan x$, $x\in (-\infty,\infty)$ has the expansion
\begin{equation}\label{equ.arcseries}
\arctan x = \sum_{n = 1}^\infty  {\frac{1}{n}\left( {\frac{x^2}{{1 + x^2 }}} \right)^{n/2} \sin \left( {n\arcsin \frac{1}{{\sqrt {1 + x^2 } }}} \right)}\,. 
\end{equation}

PROOF. Taylor's expansion for a function $f(x)$ which is analytic in an interval $I$ which includes the point $x=0$ may be written as

\begin{equation}\label{equ.nonpower}
f(x) = f(0) - \sum_{n = 1}^\infty  {\frac{{( - 1)^n }}{{n!}}x^n f^n (x)}\,. 
\end{equation}

Since $f(x)=\arctan x$ is analytic in $(-\infty,\infty)$, it has the series expansion given by~\eqref{equ.nonpower}. The derivatives of $\arctan x$ are given by~\eqref{equ.main}. The substitution of~\eqref{equ.main} in~\eqref{equ.nonpower} gives~\eqref{equ.arcseries} and the proof is complete.

\bigskip

The ratio test gives as condition for convergence of the series~\eqref{equ.arcseries}
\[
\left| {\frac{x}{{\sqrt {1 + x^2 } }}} \right| < 1\,,
\]
a condition which is automatically fulfilled for all $x$ in $(-\infty,\infty)$. 

\section{Rapidly convergent series for $\pi$ and $\pi\sqrt 3$}
In the notation of section~\ref{sec.derivative}, \eqref{equ.arcseries} can be written as
\begin{equation}\label{equ.nhds0h7}
\frac{\pi }{2} - \theta  = \sum_{n = 1}^\infty  {\frac{1}{n}} \cos ^n \theta \sin n\theta\;.
\end{equation}
A more general form of \eqref{equ.nhds0h7} can be found in~\cite{gradshtein}. What is remarkable about~\eqref{equ.nhds0h7} is that careful choices of $\theta$ yield interesting series for $\pi$ and $\pi\sqrt 3$.

\bigskip

On setting $\theta=\pi/4$, we obtain the series
\begin{equation}\label{equ.newpi1}
\frac{\pi }{4} = \sum_{n = 1}^\infty  {\frac{1}{{2^{n/2} n}}} \sin \frac{{n\pi }}{4}\,.
\end{equation}
Contrary to appearance, the right hand side contains no surd and does not require the knowledge of $\pi$ for evaluation, since $\sin( n\pi/4)$ can only take one of five possible values:
\[
\sin \frac{{n\pi }}{4} = \left\{ {\begin{array}{*{20}c}
   { - 1,n = 6,14,22,30, \ldots }  \\
   { - \frac{1}{{\sqrt 2 }},n = 5,7,13,15, \ldots }  \\
   {0,n = 4,8,12,16, \ldots }  \\
   {\frac{1}{{\sqrt 2 }},n = 1,3,9,11, \ldots }  \\
   {1,n = 2,10,18,26, \ldots }  \\
\end{array}} \right.
\]

Thus~\eqref{equ.newpi1} can be written as
\[
\frac{\pi }{4} = \sum_{n = 1}^\infty  {\frac{{( - 1)^{n - 1} }}{{4^n }}\left[ {\frac{2}{{4n - 3}} + \frac{1}{{2n - 1}} + \frac{1}{{4n - 1}}} \right]}\,, 
\]
or better still, by shifting the summation index
\begin{equation}\label{equ.gw3shg2}
\pi  = \sum_{n = 0}^\infty  {{ \frac{(-1)^n}{4^n}} \left[ {\frac{2}{{4n + 1}} + \frac{2}{{4n + 2}} + \frac{1}{{4n + 3}}} \right]} 
\,. 
\end{equation}

Formula~\eqref{equ.gw3shg2} is a base~$4$ BBP~\cite{bbp97}-type formula. The original BBP formula
\[
\pi  = \sum_{n = 0}^\infty  {\frac{1}{{16^n }}\left( {\frac{4}{{8n + 1}} - \frac{2}{{8n + 4}} - \frac{1}{{8n + 5}} - \frac{1}{{8n + 6}}} \right)}\,, 
\]
discovered using the PSLQ algorithm~\cite{ferguson99} allows the $n$th hexadecimal digit of $\pi$ to be computed without having to compute any of the previous digits and without requiring ultra high-precision~\cite{lord99}. Formula~\eqref{equ.gw3shg2} has also been obtained earlier and is listed as \mbox{formula (16)} in Bailey's compendium of known BBP-type formulas~\cite{bailey09}.

\bigskip

A converging series for $\pi\sqrt 3$ can be derived by setting $\theta=\pi/3$ in~\eqref{equ.nhds0h7}, obtaining
\[
\frac{\pi }{6} = \sum_{n = 1}^\infty  {\frac{1}{{2^n n}}\sin \frac{{n\pi }}{3}}\,. 
\]

Again since
\[
\sin \frac{{n\pi }}{3} = \frac{{\sqrt 3 }}{2} \times \left\{ {\begin{array}{*{20}c}
   {1,n = 1,2,7,8,13,14, \ldots }  \\
   {0,n = 0,3,6,9,12,15, \ldots }  \\
   { - 1,n = 4,5,10,11,16,17, \ldots }  \\
\end{array}}\,, \right.
\]
the above series can be written as
 \[
\frac{\pi }{6} = \frac{{\sqrt 3 }}{2}\sum_{n = 1}^\infty  {\frac{{( - 1)^{n - 1} }}{{2^{3n} }}\left[ {\frac{4}{{3n - 2}} + \frac{2}{{3n - 1}}} \right]}\,, 
\]
so that we have
\[
\pi\sqrt 3  = \frac{9}{8}\sum_{n = 0}^\infty  { { \frac{(-1)^n}{8^n}} \left[ {\frac{4}{{3n + 1}} + \frac{2}{{3n + 2}}} \right]}\,, 
\]
giving a base~$8$ BBP-type formula for $\pi\sqrt 3$.

\bigskip

We can obtain yet another converging series by setting $\theta=\pi/6$ in~\eqref{equ.nhds0h7}, obtaining 

\[
\frac{\pi }{3}= \sum_{n = 1}^\infty  {\left( {\frac{{\sqrt 3 }}{2}} \right)^n \frac{1}{n}\sin \frac{{n\pi }}{6}}\,,
\]
which when written out is
\begin{equation}\label{equ.newpi3}
\begin{split}
\pi\sqrt3 &= \frac{9}{64} \sum_{n = 0}^\infty  {(-1)^n\left( {\frac{3}{4}} \right)^{3n} \left[ {\frac{16}{(6n + 1)} + \frac{24}{(6n + 2)} } \right.}\\
&\qquad\qquad\qquad\qquad\left. {+\frac{24}{(6n + 3)} + \frac{18}{{6n + 4}} + \frac{9}{{6n + 5}}} \right]\,.
\end{split} 
\end{equation}
We note that technically speaking~\eqref{equ.newpi3} is not a BBP-type formula.

\section{Summary}
We have given a closed form formula for the $n^{th}$ derivative of $\arctan x$: 
\begin{align*}
\frac{{{\rm d}^n }}{{{\rm d}x^n }}\left( {\arctan x} \right) &= \frac{{( - 1)^{n - 1} (n - 1)!}}{{(1 + x^2 )^{n/2} }}\sin \left[ {n\arcsin \left( {\frac{1}{\sqrt{1 + x^2 }}} \right)} \right]\,,\nonumber\\
n&=1,2,3,\ldots
\end{align*}
We also obtained a new expansion for $\arctan x$, \mbox{$x \in ( - \infty ,\infty )$}: 
\[
\arctan x = \sum_{n = 1}^\infty  {\frac{1}{n}\left( {\frac{x^2}{{1 + x^2 }}} \right)^{n/2} \sin \left( {n\arcsin \frac{1}{{\sqrt {1 + x^2 } }}} \right)}\,. 
\]
Finally we derived rapidly convergent series for $\pi$ and $\pi\sqrt 3$:
\[
\pi  = \sum_{n = 0}^\infty  {{ \frac{(-1)^n}{4^n}} \left[ {\frac{2}{{4n + 1}} + \frac{2}{{4n + 2}} + \frac{1}{{4n + 3}}} \right]} 
\,, 
\]
\[
\pi\sqrt 3  = \frac{9}{8}\sum_{n = 0}^\infty  { { \frac{(-1)^n}{8^n}} \left[ {\frac{4}{{3n + 1}} + \frac{2}{{3n + 2}}} \right]} 
\]
and
\[
\begin{split}
\pi\sqrt3 &= \frac{9}{64} \sum_{n = 0}^\infty  {(-1)^n\left( {\frac{3}{4}} \right)^{3n} \left[ {\frac{16}{(6n + 1)} + \frac{24}{(6n + 2)} } \right.}\\
&\qquad\qquad\qquad\qquad\left. {+\frac{24}{(6n + 3)} + \frac{18}{{6n + 4}} + \frac{9}{{6n + 5}}} \right]\,.
\end{split} 
\]
The generator of the BBP-type series is the formula
\[
\frac{\pi }{2} - \theta  = \sum_{n = 1}^\infty  {\frac{1}{n}} \cos ^n \theta \sin n\theta\;.
\]

\smallskip

\textbf{Acknowledgments.}\\ 
KA is grateful to Prof. Angela Kunoth for encouragements. The authors wish to acknowledge the anonymous reviewer, whose useful comments helped to improve the quality of this work.

\smallskip 


\end{document}